\documentclass[12pt]{amsart}
\usepackage{amscd}
\newtheorem{thm}{Theorem}
\newtheorem{prop}[thm]{Proposition}
\newtheorem{lem}[thm]{Lemma}

\theoremstyle{remark}
\newtheorem{rem}[thm]{Remark}
\theoremstyle{definition}

\newtheorem{ex}[thm]{Example}

\newcommand{\C}{\mathbb{ C}}

\newcommand{\N}{\mathbb{ N}}
\newcommand{\Q}{\mathbb{ Q}}
\newcommand{\R}{\mathbb{ R}}

\newcommand{\dt}{\operatorname{det}}

\newcommand{\rk}{\operatorname{rk}}

\newcommand{\Ker}{\operatorname{Ker}}

\newcommand{\Z}{\mathbb{ Z}}
\newcommand{\cp}{{\C P^2}}

\title{On rational homotopy of four-manifolds}
\author{S.~Terzi\'c}
\address{Faculty of Science, University of Montenegro, Cetinjski put bb,
81000 Podgorica, Montenegro, Serbia and Montenegro}
\email{sterzic@cg.ac.yu}
\thanks{The paper is written while the author was postdoc supported by the
{\it DFG Graduiertenkolleg ``Mathematik im Bereich ihrer
Wechselwirkung mit der Physik''} at the Mathematical
Department of the Ludwig-Maximilians University in Munich.
The author is a member of EDGE, Research Training Network
HPRN-CT-2000-00101, supported by The European Human Potential Programme.}

\date{\today; MSC 53C25, 57R57, 58A14, 57R17}

\begin{document}

\begin{abstract}
We give explicit formulas for the ranks of the third and fourth  homotopy groups of all oriented closed simply connected four-manifolds in terms of their second Betti numbers. We also show that the rational homotopy type of these manifolds is classified by their rank and signature.
\end{abstract}

\maketitle
\section{Introduction}

In this paper we consider the problem of computation of the rational
homotopy groups and the problem of rational homotopy classification of simply
connected closed four-manifolds. Our main results could be collected as follows.
\begin{thm} \label{homotop}
Let $M$ be a closed oriented simply connected four-manifold
and $b_{2}$ its second Betti number. Then:

\begin{enumerate}
\item  If $b_{2} = 0$ then  $\rk \pi _{4}(M) = \rk \pi _{7}(M) = 1$ and
        $\pi _{p}(M)$ is finite for $p\neq 4,7$ \ ,
\item  If $b_{2} = 1$ then $\rk \pi _{2}(M) = \rk \pi _{5}(M) = 1$ and
       $\pi _{p}(M)$ is finite for $p\neq 2,5$ \ ,
\item  If $b_{2} = 2$ then $\rk \pi _{2}(M) = \rk \pi _{3}(M) = 2$ and
       $\pi _{p}(M)$ is finite for $p\neq 2,3$ \ ,
\item  If $b_{2} > 2$ then $\dim \pi _{*}(M)\otimes \Q = \infty$  and
       \[
            \rk \pi _{2}(M)=b_{2},\;\;
            \rk \pi _{3}(M)=\frac{b_{2}(b_{2}+1)}{2}-1,\;\;
            \rk \pi _{4}(M)=\frac{b_{2}(b_{2}^{2}-4)}{3} \ .
       \]
\end{enumerate}
\end{thm}

When the second Betti number is $3$, we can prove a little more.

\begin{prop}\label{small}
If $b_{2}=3$ then $\rk \pi _{5}(M)=10$.
\end{prop}

Regarding rational homotopy type classification of
simply connected closed four-manifolds, we obtain  the following.

\begin{thm}\label{type}
The rational homotopy type of a closed oriented simply connected four-manifold
is classified by its rank and signature.
\end{thm}

The  stated results are new, although as it will be clear later, we obtained
them easily by known methods. Namely, as far as we know, they can
not be found in the  well known publications presenting  the results
on topology of four-manifolds~\cite{M},~\cite{F},~\cite{DK}, nor in those presenting
the results on application of rational homotopy theory~\cite{FHT}.

\subsection{Some applications.}

\begin{rem}
The first four-manifolds we would like to apply  Theorem~\ref{homotop} are homogeneous spaces. For them the ranks of the homotopy groups are already well known. More precisely, simply connected four dimensional Riemannian homogeneous spaces  are classified in~\cite{ish} and the only such ones are $\R ^{4}$, $S^{4}$, $\C P^{2}$, $S^{2}\times S^{2}$, $\R \times S^{3}$ and $\R ^{2}\times S^{2}$.
\end{rem}

\begin{ex}
A smooth hypersurface $S_{d}$ in $\C P^{3}$ is the zero set of a homogeneous polynomial of degree $d$ in four variables. It is simply connected and $b_{2}(S_{d}) = d(6-4d+d^{2})-2$, see~\cite{schafar}. Thus, for $d\neq 1$, using Theorem~\ref{homotop}, we get
\[
\rk \pi _{2}(S_{d}) = d(6-4d+d^{2}) -2,
\]
\[
\rk \pi _{3}(S_{d}) = \frac{d(6-4d+d^{2})(d(6-4d+d^{2}) - 3)}{2} \ ,
\]
\[
\rk \pi _{4}(S_{d}) = \frac{d(6-4d+d^{2})(d^{2}(6-4d+d^{2})^{2}-6d(6-4d+d^{2})+8)}{3} \ .
\]

\end{ex}

\begin{ex}
It is known, see~\cite{schafar}, that $S_{4}$ is an example of a $K3$ surface and, thus,
\[
\rk \pi _{2}(K3) = 22, \;\; \rk \pi _{3}(K3) = 252, \;\; \rk \pi _{4}(K3) = 3520 \ .
\]

The same is also true for the ranks of homotopy groups of any logarithmic transform $L_{a}(m)L_{b}(n)(S)$ of an elliptic $K3$ surface $S$, where $m$ and $n$ are odd and relatively prime. This follows from the results of Kodaira~\cite{kod}, since he proved that such surfaces are (and only them) homotopy $K3$ surfaces.

\end{ex}

\begin{ex}
Let us consider complete intersection surfaces, i.e.,  surfaces $M$ in $\C P^{n+2}$ which are transversal intersections of $n$ hypersurfaces $Y_{1},\ldots ,Y_{n}$ that are smooth at the points of intersections. If $\deg Y_{i} = d_{i}$  then $(d_{1},\ldots ,d_{n})$ is said to be the type of $M$ and $M$ is usually denoted by   $S(d_{1},\ldots ,d_{n})$. Then, see~\cite{schafar}, $M$ is simply connected and $b_{2}(M) = e(M) -2$, where
\[
e(M) = [ {n+3\choose 2} - (n+3)\sum_{i=1}^{n}d_{i} + \sum _{i=1}^{n}d_{i}^{2} + \sum _{i\neq j} d_{i}d_{j}]\prod _{i=1}^{n}d_{i} \ .
\]
Theorem~\ref{homotop} implies that for $e(M) \geq 4$
\[
\rk \pi _{2}(M) = e-2, \;\; \rk \pi _{3}(M) = \frac{e(e-3)}{2}, \;\;
\rk \pi _{4}(M) = \frac{e(e-2)(e-4)}{3} \ .
\]

\end{ex}

\subsection{The method of the proof.}
Our proof is based on  Sullivan's
minimal model theory. It is a well  known that  Sullivan's minimal model of a simply connected space
$X$ of finite type contains complete information on the ranks of its homotopy groups,
and furthermore classifies its rational homotopy type. If the space $X$ is formal in the sense of Sullivan,
then its minimal model
coincides with the minimal model of its cohomology algebra.
By~\cite{MN} all closed oriented simply connected four-manifolds are formal,
and, thus, their cohomology algebra contains complete information on their rational homotopy.
Following this, in Section~\ref{cs} we first recall some results on real cohomology structure
of simply connected closed four-manifolds and then  in Section~\ref{main}
state the necessary background from Sullivan's minimal model theory and prove Theorem~\ref{homotop}
and Proposition~\ref{small}.
In Section~\ref{RHT} we prove  Theorem~\ref{type}.

\medskip
\noindent
{\sl Acknowledgment:} I am grateful to Yuri Petrovich Solovyov for getting me
interested into this problem. I would also like to thank Dieter Kotschick
for useful conversations.

\section{Real cohomology structure of closed oriented simply connected four-manifolds}\label{cs}

We denote by $M$  a closed oriented simply connected topological four-manifold. The symmetric bilinear form
$$
Q_{M} : H^{2}(M,\Z )\times H^{2}(M,\Z )\to \Z \ ,
$$
defined by
\[
Q_{M}(a, b ) = \langle a\cup b, [M]\rangle \
\]
is called the intersection form of $M$. Poincar\'e duality implies that (for $b_{2}\neq 0$)
the form
$Q_{M}$ is non-degenerate, and, furthermore, it is unimodular ($\dt Q_{M}=1$).
For simplicity, we will denote the cup product $a\cup b$ by $ab$ below.

The intersection form $Q_{M}$ can be diagonalised over $\R$, with $\pm 1$ as the  diagonal elements. Following standard notation (which comes from Hodge theory in the smooth case), we denote by $b_{2}^{+}$ the number of $(+1)$   and by $b_{2}^{-}$ the number of $(-1)$ in the diagonal form for $Q_{M}$. Then  $b_{2} = b_{2}^{+} + b_{2}^{-}$ and $\sigma = b_{2}^{+}-b_{2}^{-}$ are the rank and the signature of
the manifold $M$, respectively.


Using the intersection form one can easily get  an explicit description of the real cohomology algebra  of a closed oriented simply connected four-manifold.

\begin{lem}\label{b22}
Let $M$ be a closed oriented simply connected  four-manifold, such that
$b_{2}(M)\geq 2$. Then
\[
H^{*}(M,\R )\cong \R [x_{1},\ldots ,x_{b_{2}^{+}},x_{b_{2}^{+}+1},\ldots ,x_{b_{2}}]/\mathrm{relations} \ ,
\]
where $\deg x_{i} = 2$ and the $\mathrm{relations}$ are as follows:
\begin{equation}
x_{1}^{2} =\ldots = x_{b_{2}^{+}}^{2} = -x_{b_{2}^{+}+1}^{2} = \ldots = -x_{b_{2}}^{2} \ ,
\end{equation}
\begin{equation}
x_{i}x_{j} = 0, \;\; i\neq j \ .
\end{equation}
\end{lem}

\begin{proof}
We recall here the standard proof. Let $x_{i}$, $1\leq i\leq b_{2}$, be the cohomology classes  in $H^{2}(M, \R)$, representing the basis in which the intersection form $Q_{M}$ is diagonalisable. This means that
\begin{equation}\label{SD}
Q(x_{i}, x_{i}) = 1, \; \mbox{for} \; 1\leq i\leq b_{2}^{+} \ ,
\end{equation}
\begin{equation}\label{ASD}
Q(x_{i}, x_{i}) = -1, \; \mbox{for} \; b_{2}^{+}+1\leq i\leq b_{2} \ ,
\end{equation}
\begin{equation}\label{mix}
Q(x_{i}, x_{j}) = 0, \; \mbox{for} \; i\neq j \ .
\end{equation}

Denote by $V$ the generator in $H^{4}(M, \R )$ such that $\langle V, [M]\rangle = 1$. Since $x_{i}x_{j} = c \cdot V$, for some $c\in \R$, then~\eqref{SD} implies that $x_{i}^{2} = V$ for $1\leq i\leq b_{2}^{+}$, and~\eqref{ASD} implies that $x_{i}^{2} = -V$ for $b_{2}^{+}+1\leq i\leq b_{2}$. Also,~\eqref{mix} implies that $x_{i}x_{j} = 0$ for $i\neq j$.
\end{proof}

For $b_{2}=0$ or $b_{2}=1$, obviously $M$ has the following cohomology structure.

\begin{lem} \label{b2=1}
Let $M$ be a closed oriented simply connected  four-manifold.
\begin{itemize}
\item If $b_{2}(M) =1$ then
$H^{*}(M,\R ) = \R [x]$, where $\deg x=2$ and $x^{3} = 0$.
\item If $b_{2}(M) = 0$ then
$H^{*}(M,\R ) = \R [x]$, where $\deg x = 4$ and $x^{2} = 0$.
\end{itemize}
\end{lem}

\begin{rem}
As we will see in Section~\ref{RHT}, the above statements on the cohomology structure of four-manifolds are true over $\Q$ as well.
But, it will be clear from below, that for the purpose of computation of the
ranks of the homotopy groups or  determining formality, it is sufficient to work
with real coefficients.
\end{rem}

\section{The ranks of the homotopy groups of closed oriented simply connected
four-manifolds}\label{main}

\subsection{General remarks.}

We refer to~\cite{FHT} for a  comprehensive general reference for rational homotopy theory.

Let $(\mathcal{A},d_{\mathcal{A}})$ be a connected $(H^{0}(\mathcal{A}, d_{\mathcal{A}}) = k)$
and simply connected ($H^{1}(\mathcal{A}, d_{\mathcal{A}}) =0$) commutative $\N$-graded differential algebra over a field $k$ of characteristic zero. Let us consider the free $\N$-graded commutative differential algebra $(\wedge V, d)$ for a $\N$-graded vector space $V$ over $k$. We
say that $(\wedge V, d)$ is a minimal model for $(\mathcal{A},d_{\mathcal{A}})$
if $d(V)\subset \wedge ^{\geq 2}V$ and there exists a morphism
$$
f : (\wedge V, d)\to (\mathcal{A}, d_{\mathcal{A}}) \ ,
$$
which induces an isomorphism in cohomology.

Let $X$ be a simply connected topological space of  finite type.  We define the minimal model $\mu (X)$ for $X$ to be the minimal model for the
algebra $\mathcal{A}_{PL}(X)$. One says that two simply connected spaces have the same rational homotopy type if and only if there is a third space to which
they both map by maps inducing isomorphism in rational cohomology.   Then the following facts are well known. The minimal model $\mu (X)$ of a  simply connected topological space $X$ of finite type is unique up to isomorphism (which is well defined up to homotopy), it classifies the rational homotopy type of $X$ and,
furthermore, it contains complete information on the ranks of the homotopy groups of $X$. More precisely,
\begin{equation}\label{ranks}
\rk \pi _{r}(X) = \dim (\mu (X)/\mu ^{+}(X)\cdot \mu ^{+}(X))_{r}, \;
r\geq 2 \ ,
\end{equation}
where by $\mu ^{+}(X)$ we denote the elements in $\mu (X)$ of
positive degree and $\cdot$ is the usual product in $\mu (X)$.
One says that $X$ is formal in the sense of Sullivan if its minimal model
coincides with the minimal model of its cohomology algebra $(H^{*}(X,\Q ), d=0)$ (up to isomorphism).

It follows that, in the case of formal simply connected topological spaces of  finite type, we can get the ranks of their homotopy groups from their cohomology algebras, by some formal procedure. This formal procedure is, in fact, a  procedure of constructing of the  minimal model for the corresponding cohomology algebra.

\begin{rem}\label{good}
For some spaces with special  cohomology one can easily compute
their minimal models. Namely, using the terminology of~\cite{BG}, one says that
 $X$ has good cohomology if
$$
H^{*}(X,\Q )\cong \Q [x_1,\ldots ,x_n]/\langle P_1,\ldots ,P_k \rangle \ ,
$$
where the polynomials $P_1,\ldots ,P_k$ are without relations
in $\Q [x_1,\ldots ,x_n]$
(i.e. $\langle P_1,\ldots ,P_k\rangle$ is a Borel ideal).
Then in~\cite{BG} it is proved that such a space $X$ is formal and its
minimal model is given by
$$
\mu (X)=\Q [x_1,\ldots ,x_n]\otimes \wedge (y_1,\ldots ,y_k) \ ,
$$
$$
dx_i = 0, \;\; dy_i =P_i \ .
$$
Clearly,~\eqref{ranks} implies that these spaces are rationally elliptic,
i. e. $\dim \pi _{*}(X)\otimes \Q <\infty$.
\end{rem}

Unfortunately, most of four-manifolds (or precisely those with $b_2 > 2$)
do not have good cohomology.

But, the results stated in the above Remark are, in fact, consequences of  a general procedure for the construction of the minimal model for a simply connected commutative differential $\N$-graded algebra. This procedure is given by the proof of the theorem which  states the existence (and also the uniqueness up to isomorphism) of the minimal model for any such algebra. We will briefly describe this procedure here, since we are going to apply it explicitly.

\subsection{ Procedure for minimal model construction.}\label{PR}
In the procedure for the computation of the  minimal model for a simply connected commutative differential
$\N$-graded algebra $(\mathcal{A}, d)$ one starts by choosing
$\mu _{2}$ and
$m_{2} : (\mu _{2}, 0)\to (\mathcal{A}, d)$
such that
$m_{2}^{(2)} : \mu _{2}\to H^{2}(\mathcal{A}, d)$
is an isomorphism. In the inductive step, supposing that
$\mu_{k}$ and $m_{k} : (\mu _{k}, d)\to (\mathcal{A}, d)$
are constructed we extend it to $\mu _{k+1}$ and
$m_{k+1} : (\mu _{k+1}, d)\to (\mathcal{A}, d)$ with
\begin{equation}\label{MM}
\mu _{k+1}=\mu _{k}\otimes \mathcal{L}(u_{i},v_{j}) \ ,
\end{equation}
where $\mathcal{L}(u_{i},v_{j})$  denotes the vector space spanned by the elements
$u_{i}$ and $v_{j}$ corresponding to $y_{i}$ and $z_{j}$ respectively.
The latter are given by
\begin{equation}\label{M1}
H^{k+1}(\mathcal{A})=\Im m_{k}^{(k+1)}\oplus \mathcal{L}(y_{i})
\end{equation}
and
\begin{equation}\label{M2}
\Ker m_{k}^{(k+2)}=\mathcal{L}(z_{j}) \ .
\end{equation}
Then we have that $m_{k}(z_{j})=dw_{j}$ for some $w_{j}\in \mathcal{A}$ and the homomorphism $m_{k+1}$ is defined by
$m_{k+1}(u_{i})=y_{i}$, $m_{k+1}(v_{j})=w_{j}$ and $du_{i}=0$, $dv_{j}=z_{j}$.

\begin{rem}
In general, for a simply connected topological space $X$ we have that  $ \mathcal{A} = \mathcal{A}_{PL}(X)$ and, obviously, by~\eqref{ranks}, we see that   $\rk \pi _{k+1}(X)$ is the number of generators in the above procedure we add to $\mu _{k}(X)$, in order to obtain $\mu _{k+1}(X)$ .
\end{rem}

Note that for a formal $X$, the algebra  $\mu (X)\otimes _{\Q}k$ coincides with the minimal model of the cohomology algebra $(H^{*}(X,k), d=0)$ for any field $k$ of characteristic zero. The converse is also true. If there exists a field $k$ of characteristic zero for which $\mu (X)\otimes _{\Q}k$ is the minimal model for the cohomology algebra $(H^{*}(X,k), d=0)$, then $X$ is formal.

\begin{rem}
Obviously,~\eqref{ranks} implies that for the purpose of calculating the ranks of the homotopy groups of $X$ we can use $\mu (X)\otimes _{\Q}\R$ as well.
In the case of formal $X$ it means that we can apply the above procedure to
$H^{*}(X,\R )$.
\end{rem}

\subsection{Computation of the ranks of the homotopy
groups.}\label{computation}

Before we proceed to the computation of the ranks of the low degree homotopy groups of four-manifolds,
let us note the following important facts.

\begin{rem}\label{hyp}
All closed oriented simply connected four-manifolds with $b_{2}> 2$ are rationally hyperbolic, i.e.
$\dim \pi _{*}(M)\otimes \Q =\infty$. One can see that using the fact that for $M$ rationally
elliptic must be satisfied $\sum 2k\rk \pi _{2k}(M)\leq \dim M$ (see~\cite{FHT}). Then from Hurewicz
isomorphism it follows that $2b_{2}\leq \dim M$. In particular, it implies that for $b_{2}> 2$
these spaces do not have good cohomology.
 \end{rem}

\begin{rem}\label{form}
All closed oriented simply connected four-manifolds are formal in the sense of Sullivan. One can see it using the results of~\cite{MN} which say that any
compact simply connected manifold of dimension $\leq 6$ is formal.
\end{rem}

The  Remarks~\ref{hyp} and~\ref{form} together with  Procedure~\ref{PR}  for minimal model computation and
the knowledge of cohomology structure of four-manifolds, make us possible to prove the Theorem~\ref{homotop}.

\begin{rem}
Note that, using Hurewicz isomorphism, we already know that $\rk \pi _{2}(M) = b_{2}$.
\end{rem}

{\bf Proof of the Theorem~\ref{homotop}.}
Because of formality the minimal model of a closed oriented simply connected four-manifold $M$ is the
minimal model of its cohomology algebra.
Therefore, to compute the minimal model $\mu (M)\otimes _{\Q}\R$, we can
apply Procedure~\ref{PR} to the algebra $(\mathcal{A}, d) = (H^{*}(M,\R ),
d=0)$. For simplicity we denote $\mu (M)\otimes _{\Q}\R$  by $\mu (M)$
and $H^{*}(M,\R )$ by $H^{*}(M)$ below.
As stated in the theorem we distinguish the following cases.
\\
{\bf 1.} For $b_{2} = 0$ Lemma~\ref{b2=1} immediately implies that
       $\mu (M) = \wedge    (x, u)$, where $\deg x = 4$, $\deg u = 7$ and
       $dx = 0$, $du = x^{2}$. Thus,
       $$
        \rk \pi _{4}(M) = \rk \pi _{7}(M) = 1
       $$
       and $\pi _{p}(M)$ are finite for $p\neq 4,7$.\\
{\bf 2.} For $b_{2}= 1$ it follows from Lemma~\ref{b2=1} that
         $\mu (M) = \wedge (x, u)$, where $\deg x = 2$, $\deg u = 5$ and
       $dx = 0$,  $du = x^{3}$. Thus,
       $$
       \rk \pi _{2}(M) = \rk \pi _{5}(M) = 1
       $$
       and $\pi _{p}(M)$ are finite for $p\neq 2, 5$.\\
{\bf 3.} Let the second Betti number of $M$ be $2$. Lemma~\ref{cs} implies
       that $M$ has good cohomology, since the polynomials $x_1^2\pm x_2^2$,
       $x_1x_2$ are without relations in $\R [x_1,x_2]$. By Remark~\ref{good},        the minimal model for $M$ is given by
       $$
       \mu (M) = \R [x_1,x_2]\otimes \wedge (y_1,y_2), \; dx_i=0,\;
       dy_{1}=x_{1}^{2}\pm x_{2}^{2},\; dy_{2}=x_{1}x_{2} \ .
       $$
       Thus, $\pi _{p}(M)$ are finite for $p\neq 2,3$ and
       $$
         \rk \pi _{2}(M)=\rk \pi _{3}(M) = 2 \ .
       $$
{\bf 4.} Let $b_{2} > 2$. Remark~\ref{hyp} implies that in this case
       $\dim \pi _{*}(M)\otimes \Q = \infty$.
       We will use here the results on cohomology of $M$ proved in
       Lemma~\ref{b22}.
       According to Procedure~\ref{PR} for minimal model construction,
       it follows that
       $\mu _{2}(M)=\R [x_1,\ldots ,x_{b_2}]$, $m_{2}(x_{i})=[x_{i}]$.
       Therefore,
       $\rk \pi _{2}(M) = b_{2}$.
At the next step in the application of Procedure~\ref{PR}, we know that $H^{3}(M)=0$ and
$\Ker m_{2}^{(4)}=\mathcal{L}(x_{1}^{2}\pm x_{i}^{2}, x_{i}x_{j})$,
$2\leq i\leq b_{2}$, $1\leq i< j\leq b_{2}$.
We take here $b_{2}^{+}-1$ times sign $(-)$ and $b_{2}^{-}$ times sign $(+)$.
Since the elements $x_{1}^{2}\pm x_{i}^{2}$, $x_{i}x_{j}$ are linearly independent, we obtain that
$$
\mu _{3}(M)=\mu _{2}(M)\otimes \mathcal{L}(v_{i},v_{ij}), \;
2\leq i\leq b_{2}, 1\leq i< j\leq b_{2} \ ,
$$
$$
 dv_{i}=x_{1}^{2}\pm x_{i}^{2}, \; dv_{ij}=x_{i}x_{j}, \; m_{3}(v_{i})=m_{3}(v_{ij})=0 \ .
$$
This  implies that
\[
\rk \pi _{3}(M)= \frac{b_{2}(b_{2}+1)}{2} -1 \ .
\]

In order to continue this procedure, let us note the following. For $k\geq 3$, $\mu _{k+1}(M)$ is given by
$\mu _{k+1}(M)=\mu _{k}(M)\otimes \mathcal{L}(w_j)$, where $w_j$ correspond to basis for $H^{k+2}(\mu _{k}(M), d)$. We get this from~\eqref{MM} using the following two observation.

First,~\eqref{M2} implies that $\Ker m_{k}^{(k+2)}(M)=H^{k+2}(\mu _{k}(M), d)$ for $k\geq 3$, since then $H^{k+2}(M)=0$. Second, since
$\Im m_{3}^{(4)}\cong \Im m_{2}^{(4)}/\Ker m_{2}^{(4)}\cong H^{4}(M)$, it follows that, for $k\geq 3$, there are no $y_{i}$'s in~\eqref{M1}.

Thus, in order to construct $\mu _{4}(M)$, we need to find the basis for $H^{5}(\mu _{3}(M))$.
Since in $\mu _{3}(M)$ we have no nontrivial 5 dimensional coboundaries, this is equivalent to finding
the basis for the 5 dimensional cocycles in $\mu _{3}(M)$.
 Any cochain of degree $5$ in $\mu _{3}(M)$ is of the form
$$
c=\sum_{i=2}^{b_2}P_{i}v_{i}+\sum_{1\leq i<j\leq b_2}P_{ij}v_{ij} \ ,
$$
where $P_{i}=\sum_{k=1}^{b_2}\alpha _{i}^{k}x_{k}$ and $P_{ij}=\sum_{k=1}^{b_2}\alpha _{ij}^{k}x_{k}$.
Computing the coefficients for $x_{1}^{3}$, $x_{i}^{3}$, $2\leq i\leq b_2$,
$x_{1}^{2}x_{j}$, $2\leq j\leq b_2$, $x_{i}^{2}x_{j}$, $i\neq j$, $2\leq i\leq b_2$, $1\leq j\leq b_2$,
 $x_{i}x_{j}x_{k}$, $1\leq i< j< k\leq b_2$ in the expression for
 $d(c)$, one obtains that the equation $d(c)=0$ gives rise to the
 following system of linear equations (respectively).
\begin{equation}\label{jedan}
        \sum _{i=2}^{b_2}\alpha _{i}^{1}=0,
\end{equation}
\begin{equation}\label{dva}
         \alpha _{j}^{j}=0, \; 2\leq j\leq b_2 \ ,
 \end{equation}
\begin{equation}\label{tri}
         \sum _{i=2}^{b_{2}}\alpha _{i}^{j} + \alpha _{1j}^{1} = 0,\;
          2\leq j\leq b_2 \ ,
\end{equation}
\begin{equation}\label{cetiri}
-\alpha _{i}^{j}+\alpha _{ij}^{i}=0,\; i<j, \; 2\leq i\leq b_{2}^{+},
3\leq j\leq b_{2} \ ,
\end{equation}
$$
 -\alpha _{i}^{j}+\alpha _{ji}^{i}=0, \; i>j, \; 2\leq i\leq b_{2}^{+}, 1\leq j\leq b_{2} -1 \ ,
$$
$$
 \alpha _{i}^{j}+\alpha _{ij}^{i}=0,\; i<j, \; b_{2}^{+}+1\leq i\leq b_{2}, b_{2}^{+} + 2\leq j\leq b_{2} \ ,
$$
$$
  \alpha _{i}^{j}+\alpha _{ji}^{i}=0, \; i>j, \; b_{2}^{+}+1\leq i\leq b_{2}, 1\leq j\leq b_{2}-1 \ ,
$$
\begin{equation}\label{pet}
\alpha _{ij}^{k}+\alpha _{ik}^{j}+\alpha _{jk}^{i}=0,\; 1\leq i< j< k\leq b_{2} \ .
\end{equation}

The number of variables in the above system is
$\frac{b_{2}(b_{2}-1)(b_2+2)}{2}$. By the inspection one sees that
the equation~\eqref{jedan} eliminates $1$ variable, each of the systems~\eqref{dva} and~\eqref{tri}
eliminates $b_{2}-1$ variables, the system~\eqref{cetiri} eliminates  $(b_2-2)(b_2-1)+b_2-1$ variables
and the system~\eqref{pet} eliminates $b_{2}\choose 3$ variables.
Thus, the dimension of
the solution space for the above system is $\frac{b_{2}(b_{2}-1)(b_2+2)}{2} -
(b_{2}^{2} + {b_{2}\choose 3}) = \frac{b_2(b_2^{2}-4)}{3}$.
So, at this  step in the construction of the minimal model,
we extend $\mu _{3}(M)$ by adding the  generators
$w_{k}$, $1\leq k\leq  \frac{b_2(b_2^{2}-4)}{3}$ of degree $4$, i. e.
$$
\mu _{4}(M)=\mu _{3}(M)\otimes \mathcal{L}(w_{k}), \;
1\leq k\leq  \frac{b_2(b_2^{2}-4)}{3} \ .
$$
The differentials $dw_{k}$ are given by some basis
for the solution space
of the above system. Thus, we have
\[
\rk \pi _{4}(M)=\frac{b_2(b_2^{2}-4)}{3} \ .
\]

One can continue the above procedure for the  construction of the minimal model and calculation the ranks of the homotopy groups, but it is obvious that at
each  step the vector space $H^{k+2}(\mu _{k}(M), d)$ for which we want to
get a  basis becomes bigger and more complicated.
Clearly,  $\rk \pi _{k+1}(M)$ is given by the dimension of the $H^{k+2}(\mu _{k}(M), d)$ and it is some polynomial $P_{k+1}(b_{2})$ in $b_{2}(M)$.
In order to continue the process we need to have a basis for
$H^{k+2}(\mu _{k}(M), d)$ as well.

\begin{rem}
On each step in construction of the minimal model we should solve some system of linear
equations, whose dimension of the solution space determines the number of generators of corresponding
degree in the minimal model and the differentials for these generators are given by some solution of that
system. The calculation procedure done in the proof of Theorem~\ref{homotop}
suggests that the dimension of the solution space for such system does not depend on the
signature of the manifold, but only on its rank, while its explicit solution does.
In other words,  it suggests that the ranks of the homotopy groups of simply connected four-manifold
are completely determined by the rank of its intersection form.
Our attempt to obtain explicit proof for it
following the proof of Theorem~\ref{homotop} involved  complicated
calculations, which we were not able to carry out.
\end{rem}

\begin{rem}
Note that for a  simply connected topological space $X$ of finite type and
finite rational Lusternik-Schnirelman category, hyperbolicity of $X$ implies
that its rational homotopy groups grow exponentially ( so called rational dichotomy,~\cite{FH}). This explains why one should expect the computations in the
 above procedure to be more and more complicated.
This also gives that we can not expect to control the degrees of the
polynomials $P_{k}(b_{2})$ with the growth of $k$.
\end{rem}

{\bf Proof of the Proposition~\ref{small}.}
If we continue the procedure we started with in the  proof of Theorem~\ref{homotop},
we need to extend $\mu _{4}(M)$ by adding generators of degree $5$ that correspond
to basis of $H^{6}(\mu _{4}(M))$. First note that any cocycle of degree $6$ in $\mu _{4}(M)$ is cohomologous to the cocycle of the form
\begin{equation}\label{ran}
c=\sum \alpha _{ij}v_{i}v_{j}+\sum \alpha _{ijk}v_{i}v_{jk}+\sum \alpha _{ijkl}v_{ij}v_{kl}+\sum \beta _{ij}x_{i}w_{j} \ .
\end{equation}
To simplify the calculations we can assume that $b_{2}^{+}$ is also 3,
since it  is clear that the dimension of the solution space for the equation
$d(c)=0$ does  not depend on $b_{2}^{+}$. According to the proof of the above theorem we can
take the differentials $dw_i$ to be as follows.
$$
dw_1 = x_1v_2 - x_1v_3 + x_2v_{12} - x_3v_{13} \ ,
$$
$$
dw_2 = x_3v_2 - x_1v_{13} + x_2v_{23} \ ,
$$
$$
dw_3 = x_2v_3 - x_1v_{12} + x_3v_{23} \ ,
$$
$$
dw_4 = x_3v_{12} - x_1v_{23}, \;\; dw_5 = x_2v_{13} - x_1v_{23} \ .
$$

On that way, for $c$ of the form~\eqref{ran}, the differential $d(c)$ is given
as follows.
\begin{eqnarray}
d(c) & = & (-\alpha _{23} + \beta _{11})x_{1}^{2}v_2 +
(\alpha _{23} + \beta _{32})x_{3}^{2}v_2 + (-\alpha _{212} + \beta _{21})x_1x_2v_2
 + \nonumber \\
& &  (-\alpha _{223} + \beta _{22})x_2x_3v_2 + (-\alpha _{213} + \beta _{31} + \beta _{12})x_1x_3v_2 +
\nonumber \\
& & (\alpha _{23} - \beta _{11})x _{1}^{2}v_3 + (-\alpha _{312} - \beta _{21} + \beta
_{13})x_1x_2v_3 +
\nonumber \\
& &  (-\alpha _{23} + \beta _{23}) x_{2}^{2}v_3 + (\alpha _{212} + \alpha _{312} -
\beta _{13})x_{1}^{2}v_{12} + \nonumber \\
& &
(-\alpha _{313} - \beta _{31})x_1x_3v_3 + (-\alpha _{1213} - \beta _{33} + \beta
_{14})x_1x_3v_{12} + \nonumber \\
& & (-\alpha _{323} + \beta _{33})x_2x_3v_3 +  (-\alpha _{1223} + \beta _{31} + \beta
_{24})x_2x_3v_{12} + \nonumber \\
& &  (-\alpha _{212} + \beta _{21})x_{2}^{2}v_{12} + (\alpha _{213} +\alpha _{313} -\beta
_{12})x_{1}^{2}v_{13} + \nonumber \\
& & (-\alpha _{312} + \beta_{34})x_{3}^{2}v_{12} +
 (\alpha _{1213} - \beta _{22} + \beta
_{15})x_1x_2v_{13} + \nonumber \\
& &
  (-\alpha _{213} + \beta
_{25})x_{2}^{2}v_{13} + (-\alpha _{1323} - \beta _{21} + \beta
_{35})x_2x_3v_{13} + \nonumber \\
& & (-\alpha _{313} -\beta
_{31})x_{3}^{2}v_{13} +  (\alpha _{223} + \alpha _{323} - \beta _{14}
- \beta _{15})x_{1}^{2}v_{23} + \nonumber \\
& &
(-\alpha _{223} + \beta
_{22})x_{2}^{2}v_{23} +  (\alpha _{1223} +\beta _{12} - \beta
_{24} - \beta _{25})x_1x_2v_{23} + \nonumber \\
& & (-\alpha _{323} +\beta
_{33})x_{3}^{2}v_{23} +
(\alpha _{1323} + \beta
_{13}-\beta _{34} - \beta _{35})x_1x_3v_{23} + \nonumber \\
& & (\beta _{11} - \beta _{23})x_1x_2v_{12} + (-\beta _{11} -\beta
_{32})x_1x_3v_{13} + (\beta _{23} + \beta _{32})x_2x_3v_{23} \ .
\end{eqnarray}
Using the  above expression one can easily see that all $\alpha$'s can be expressed
in terms of $\beta$'s. Besides that we see that $\beta _{23} = \beta _{11}$, $\beta _{32}= -\beta
_{11}$, $\beta _{34}=\beta _{13}-\beta _{21}$, $\beta _{15}=\beta
_{22}+\beta _{33}-\beta _{14}$ and $\beta _{25} = \beta
_{31}+\beta _{12}$, while the other $\beta$'s are linearly
independent. It implies that for $b_{2} =3$ we have that  $\rk \pi _{5}(M) = 10$.

\section{On homotopy type classification}\label{RHT}

{\bf Proof of the Theorem~\ref{type}.} As we already mentioned, Sullivan's minimal model
theory provides a bijection between rational homotopy types of simply connected spaces of  finite type and
 isomorphism classes of minimal Sullivan algebras over $\Q$.
Thus, in our case, two closed oriented simply connected  four-manifolds have the same
rational homotopy type if and only if they have isomorphic minimal models over $\Q$.  Any such four-manifold
is formal and its rational cohomology structure is determined by its intersection form over  $\Q$.
It follows that the rational homotopy type of a closed oriented simply
connected four-manifold is classified  by its intersection form over $\Q$.
In other words, two closed oriented simply connected four-manifolds have the same rational
homotopy type if and only if their intersection forms are equivalent over $\Q$. Using the Hasse-Minkowski theorem
 one can see that any quadratic form over $\Q$
which has integral unimodular lattice is equivalent (over $\Q$) to some diagonal form with $\pm 1$ diagonal elements~\cite{OM}. It follows that two intersection forms are equivalent over $\Q$ if and only if they have the same rank and signature.

By a result of Pontryagin-Wall (see~\cite{M}), the homotopy type of a simply connected closed oriented
four-manifold is determined by its intersection form.

Thus, for closed oriented simply connected four-manifolds, we have the following homotopy type classification.
\begin{enumerate}
\item Rational homotopy type $\sim$ rank and signature\ .
\item Homotopy type $\sim$ intersection form\ .
\end{enumerate}

Note that not every homotopy type of a closed oriented simply connected
four-manifold can be realised by a such smooth manifold.
Namely, by a result of Freedman~\cite{F}, for any unimodular symmetric bilinear form
$Q$ there exists a closed oriented simply connected
four-manifold having $Q$ as its intersection form.
On the other hand the theorems of Rokhlin~\cite{R} and Donaldson~\cite{D}
give the constraints on the intersection form of smooth four-manifold.
This implies the existence of the intersection forms (like $E_8$)
which can not be realised by some smooth four-manifolds.

We see that, in contrast to homotopy type, every rational homotopy type
of closed oriented simply connected four-manifolds has a smooth representative.
More precisely, any closed oriented simply connected four-manifold is
rationally homotopy
equivalent to a connected sum of $\cp$ 's  and $\overline{\cp}$ 's.

\begin{rem}
One can define $\R$-homotopy type of a simply connected space $X$ of a finite
type to be the equivalence class of the algebra $\mu _{\R}(X) = \mu (X) \otimes _{\Q}\R$, where $\mu (X)$  is the minimal model for $X$, see~\cite{L}. Then obviously we have that two spaces have the same $\R$-homotopy type if and only if their minimal models are equivalent over $\R$.
Theorem~\ref{type} implies that in the class of closed oriented
simply connected four-manifolds there is no difference between the rational and the real
homotopy types.
\end{rem}

\begin{rem}
Note that, in general, in the class of simply connected spaces of finite type
we have strict inclusion among the spaces having the same rational and the same real homotopy type. To see that it is enough to construct two minimal algebras starting from two rational quadratic forms which have the same rank and signature, but which are not equivalent over the rationals. Since for any minimal algebra  $\mu $ over $\Q$ there exists  the simply connected space of  finite type having  $\mu$ as its minimal model, we get on this way two simply connected spaces of finite type which have the same real but different rational homotopy type.
\end{rem}


\bibliographystyle{amsplain}

\end{document}